\documentclass[leqno]{article}
\usepackage{amsmath,amsfonts,amsthm,amssymb,indentfirst}

\newcommand{\<}{\kern.0833em}

\newtheorem{theorem}{Theorem}
\newtheorem{lemma}[theorem]{Lemma}

\newtheorem{definition}[theorem]{Definition}
\newtheorem{question}[theorem]{Question}

\newcommand{\Z}{\mathbb Z}
\newcommand{\Q}{\mathbb Q}
\newcommand{\Sym}{\mathrm{Sym}}

\makeatletter
\newcommand{\xlabel}{\stepcounter{equation}
  \gdef\@currentlabel{\p@equation\theequation}{\rm(\@currentlabel)}}
\makeatother
\newenvironment{xlist}
  {\begin{list}{\xlabel}
    {\setlength{\rightmargin}{20pt}
     \setlength{\leftmargin}{37pt}
     \setlength{\labelsep}{20pt}
     \setlength{\labelwidth}{20pt}}}
  {\end{list}}

\begin{document}

\title{Generating infinite symmetric groups%
\thanks{
\protect\\
Preprint versions of this paper: 
http://math.berkeley.edu/\protect\linebreak[0]%
{$\!\sim$}gbergman\protect\linebreak[0]%
/papers/Sym\<\_\,Omega:1.\{tex,dvi\}
and arXiv:math.GR/0401304\,.
}}
\author{George M. Bergman}
\maketitle

\begin{abstract}
Let $S=\Sym(\Omega)$ be the group of all permutations of an infinite
set $\Omega\<.$
Extending an argument of Macpherson and Neumann, it is shown that
if $U$ is a generating set for $S$ as a group, respectively
as a monoid, then there exists a positive integer $n$ such that
every element of $S$ may be written as a group word, respectively a
monoid word, of length $\leq n$ in the elements of $U.$

Some related questions and recent results by others are noted,
and a brief proof is given of a result of Ore's on commutators
that is used in the proof of the above result.
\end{abstract}

\section{Introduction, notation, and some
lemmas on full moieties.}\label{S.moieties}

In \cite[Theorem~1.1]{DM&PN} Macpherson and Neumann show that if
$\Omega$ is an infinite set, then the group $S=\Sym(\Omega)$ is not
the union of a chain of $\leq|\<\Omega\<|$ proper subgroups.
We will repeat the beautiful proof of that result, with modifications
that will allow us to obtain along with it the result stated in the
abstract.
The present section is devoted to obtaining strengthened versions of
the lemmas used in that proof.

Following the notation of~\cite{DM&PN},
for $\Omega$ an infinite set, $\Sym(\Omega),$ generally
abbreviated $S,$ will denote the group of all permutations of $\Omega,$
and such permutations will be written to the right of their arguments.
For subsets $\Sigma\subseteq\Omega$ and $U\subseteq S,$
the symbol $U_{(\Sigma)}$ will denote
the set of elements of $U$ that stabilize $\Sigma$ pointwise, and
$U_{\{\Sigma\}}$ the set $\{f\in U:\Sigma f=\Sigma\}.$
(In~\cite{DM&PN} this notation was used only for $U$ a subgroup.)
A subset $\Sigma\subseteq\Omega$ will be called {\em full} with respect
to $U\subseteq S$ if the set of permutations of $\Sigma$
induced by members of $U_{\{\Sigma\}}$ is all of $\Sym(\Sigma).$
The cardinality of a set $X$ will be written $|X|,$
and a subset $\Sigma\subseteq\Omega$ will be called
a {\em moiety} if $|\Sigma|=|\<\Omega\<|=|\<\Omega-\Sigma\<|.$
\vspace{6pt}

Suppose $\Sigma_1$ and $\Sigma_2$ are moieties of $\Omega$ whose
intersection is also a moiety, and whose union is all of $\Omega\<.$
Then \cite[Lemma~2.3]{DM&PN} says that if $G$ is a subgroup of
$S=\Sym(\Omega)$ such that
$\Sigma_1$ and $\Sigma_2$ are both full with respect to $G,$ then $G=S.$
To strengthen this result, we will consider subsets
$U, V\subseteq S,$ closed under inverses, such that $\Sigma_1$ is full
with respect to $U$ and $\Sigma_2$ with respect to $V.$
By the lemma cited, $\langle\<U\cup V\rangle=S;$
our version of this result will bound the number of
factors from $U$ and $V$ needed to get the general element of $S.$

Our proof will use the following fact,
first proved by Ore~\cite{OO}.
Much stronger results have been proved since.
In~\S\ref{S.rplt} we will give
a self-contained proof of a statement of intermediate strength.

\begin{lemma}[{\rm \cite{OO}, cf.~\S\ref{S.rplt}
below}\textbf{}]\label{L.cmttr}
For $\Omega$ an infinite set, every element $f\in\Sym(\Omega)$ can be
written as a commutator, $f=g^{-1}h^{-1}g\<h$
$(g,h\in\Sym(\Omega)).$\qed
\end{lemma}

Here now is the result on full moieties.

\begin{lemma}[{\rm cf.~\cite[Lemma~2.3]{DM&PN}}\textbf{}]\label{L.S1S2}
Suppose $\Sigma_1,\,\Sigma_2$ are moieties of $\Omega$ whose
intersection is a moiety, and whose union is all of $\Omega\<;$ and
suppose $U,\,V$ are subsets of $S=\Sym(\Omega),$ each closed under
inverses, such that $\Sigma_1$ is full with respect to $U$ and
$\Sigma_2$ is full with respect to $V.$
Then $S=(UV)^4\<V\,\cup\,(VU)^4\<U.$
\end{lemma}\begin{proof}
Note that $S_{(\Omega\<-(\Sigma_1\<\cap\<\Sigma_2))}\!\cong
\Sym(\Sigma_1\cap\Sigma_2).$
By Lemma~\ref{L.cmttr}, any element $f$ of the latter group may be
written as a commutator $f=g^{-1}h^{-1}g\<h.$
Since $\Sigma_1$ is full with respect to $U,$ we can find
an element of $U_{\{\Sigma_1\}}$ which behaves like $g$
on $\Sigma_1\cap\Sigma_2$ and as the identity on $\Sigma_1-\Sigma_2;$
likewise, we can find
an element of $V_{\{\Sigma_2\}}$ which behaves like $h$
on $\Sigma_1\cap\Sigma_2$ and as the identity on $\Sigma_2-\Sigma_1.$
Clearly the commutator of these elements behaves like $f$ on
$\Sigma_1\cap\Sigma_2$ and as the identity
on $\Omega-(\Sigma_1\cap\Sigma_2).$
Hence

\begin{xlist}\item\label{x.UVUV}
$S_{(\Omega\<-(\Sigma_1\<\cap\<\Sigma_2))}\,\subseteq\,UVUV.$
\end{xlist}

Now $|\Sigma_1\cap\Sigma_2|=|\<\Omega\<|= |\Sigma_2-\Sigma_1|,$ where
the first equality holds because $\Sigma_1\cap\Sigma_2$ is a moiety
and the second because $\Sigma_1=\Omega-(\Sigma_2-\Sigma_1)$ is one.
Hence $\Sym(\Sigma_2)$ contains an element interchanging the
subsets $\Sigma_1\cap\Sigma_2$ and $\Sigma_2-\Sigma_1;$ hence $V$
has an element which behaves that way on $\Sigma_2$ (and in
an unspecified manner on $\Sigma_1-\Sigma_2).$
Conjugating~(\ref{x.UVUV}) by such an element, we get

\begin{xlist}\item\label{x.VUVUVV}
$S_{(\Sigma_1)}\,\subseteq\,VUVUVV.$
\end{xlist}

Since the assumptions on $\Sigma_1$ and $\Sigma_2$ are
symmetric, we also have the corresponding formula for $S_{(\Sigma_2)},$
with $U$ and $V$ interchanged.

Now suppose we are given $f\in S,$ which we wish to write as a product
of elements of $U$ and $V.$

We shall see, roughly, that a product of one element from $U$ and
one element of $V$ suffices to distribute the elements of $\Omega$
between $\Sigma_1$ and its complement exactly as $f$ does.
An application of an element of $U$ will then put the elements
that have been moved into $\Sigma_1$ in exactly the desired places, and
a final application of~(\ref{x.VUVUVV}) will administer the
{\em coup de gr\^ace}.

The details:  Note first that the
set $(\Sigma_1\cap\Sigma_2)\<f^{-1}$ must either contain $|\<\Omega\<|$
elements of $\Sigma_1$ or $|\<\Omega\<|$ elements of $\Sigma_2;$
without loss of generality assume the former.
(This is the reason for the word ``roughly'' in the preceding paragraph.
In the contrary case, the roles stated there for $\Sigma_1$
and $\Sigma_2,$ and likewise for $U$ and $V,$ will be reversed.)
In particular,
$\Sigma_1f^{-1}$ contains $|\<\Omega\<|$ elements of $\Sigma_1.$
Hence we can find a permutation
$a\in U_{\{\Sigma_1\}}$ which maps all elements of $\Sigma_1$
which are not in $\Sigma_1 f^{-1}$ (if any) into $\Sigma_1\cap\Sigma_2,$
and which also maps into that set $|\<\Omega\<|$ elements
of $\Sigma_1$ which {\em are} in $\Sigma_1f^{-1}.$
These conditions, and the fact that $a\in U_{\{\Sigma_1\}}$ takes
$\Omega-\Sigma_1$ to itself, together imply that $a$
maps all elements of $(\Omega-\Sigma_1)f^{-1}$
(both those in $\Sigma_1$ and those in $\Omega-\Sigma_1)$
into $\Sigma_2,$ and also takes $|\<\Omega\<|$
elements of $\Sigma_1 f^{-1}$ there.
We can now choose $b\in V_{\{\Sigma_2\}}$ which maps
into $\Omega-\Sigma_1$ the images under $a$ of all elements
of $(\Omega-\Sigma_1)f^{-1}$ and nothing else; i.e., such
that $(\Omega-\Sigma_1)f^{-1}ab=\Omega-\Sigma_1.$
Taking complements, we have $\Sigma_1 f^{-1} ab = \Sigma_1,$
so as $\Sigma_1$ is full with respect to $U,$
we can find $c\in U_{\{\Sigma_1\}}$ which agrees on $\Sigma_1$
with the inverse of $f^{-1} ab;$
i.e., such that $f^{-1}abc\in S_{(\Sigma_1)}.$
Now~(\ref{x.VUVUVV}) applied to the inverse of the latter element
gives us $(abc)^{-1}f\in VUVUVV,$ so $f\in (UVU)(VUVUVV)=(UV)^4 V.$
As noted earlier, the roles of $U$ and $V$ may be the opposite of those
we have assumed, giving the alternative possibility $f\in (VU)^4\<U.$
\end{proof}

The result from~\cite{DM&PN} that we have just strengthened was
used there to show that if a subgroup $G\leq S$ has a full moiety,
then there exists $x\in S$ such that $\langle G\cup\{x\}\rangle = S.$
The version proved above yields the more precise statement:

\begin{lemma}[{\rm cf.\ \cite[Lemma~2.4]{DM&PN}}\textbf{}]\label{L.G&x}
If a subset $U\subseteq S$ closed under inverses has a full moiety,
then there exists $x\in S$ of order $2$ such that
$(Ux)^7\<U^2x\,\cup\,(x\<U)^7 x\<U^2\<=\<S.$
\end{lemma}\begin{proof}
Given a full moiety $\Sigma_1$ for $U,$ choose any moiety
$\Sigma_2\subseteq\Omega$
such that $\Sigma_1\cap\Sigma_2$ is a moiety and
$\Sigma_1\cup\Sigma_2=\Omega\<.$
Since $\Omega\,{-}\,\Sigma_1$ and $\Omega\,{-}\,\Sigma_2$ are disjoint
and both have the cardinality of $\Omega,$ we can find an element $x$
of order $2$ which interchanges those two sets, and hence also
interchanges their complements, $\Sigma_1$ and $\Sigma_2.$
The fact that $\Sigma_1$ is a full moiety for $U$ makes
$\Sigma_2=\Sigma_1\<x$ a full moiety for $x^{-1}Ux=x\,Ux.$
Setting $x\,Ux=V,$ we may apply the preceding lemma.
The expression $(UV)^4V$ becomes
$(Ux\<Ux)^4 x\<Ux=(Ux)^8 x\<Ux=(Ux)^7Uxx\<Ux=(Ux)^7U^2x,$ while the
other term is the conjugate of this by $x,$ namely $(x\<U)^7 x\<U^2.$
\end{proof}

We conclude this section with a diagonal argument, using nothing but
the definition of full moiety and basic set theory, which we extract
virtually unchanged from the proof of~\cite[Theorem~1.1]{DM&PN}:

\begin{lemma}\label{L.oEfull}
Let $\Omega$ be an infinite set, let $S=\Sym(\Omega),$ and let
$(U_i)_{i\in I}$ be any family of subsets of $S$ such that
$\bigcup_I\,U_i = S$ and $|I|\leq|\<\Omega\<|.$
Then $\Omega$ contains a full moiety with respect to
at least one of the $U_i.$
\end{lemma}\begin{proof}
Since $|\<\Omega\<|$ is infinite and $I\leq|\<\Omega\<|,$ we can write
$\Omega$ as a union of disjoint moieties $\Sigma_i$ $(i\in I).$
If there are no full moieties with respect to $U_i$ for any $i,$ then
in particular, for each $i$ the set $\Sigma_i$ is non-full with respect
to $U_i,$ so we can choose $f_i\in\Sym(\Sigma_i)$ which
is not the restriction to $\Sigma_i$ of a member of $U_i.$
Now let $f\in\Sym(\Omega)$ be the permutation whose
restriction to each $\Sigma_i$ is $f_i.$
Then $f$ cannot belong to any of the $U_i,$ contradicting
the assumption $\bigcup_I\,U_i = S,$ and completing the proof.
\end{proof}

\section{Chains of subsets of $\Sym(\Omega).$}\label{S.U^n}

Let us begin by recovering~\cite[Theorem~1.1]{DM&PN}.
Our statement will be the contrapositive of that in~\cite{DM&PN}.

I also include parenthetically the corresponding statement with
chains of submonoids in place of chains of subgroups.
As we will see, this follows
trivially from the result on chains of subgroups;
but it took me a long time to discover that trivial argument.

\begin{theorem}[{\rm\cite[Theorem~1.1]{DM&PN}}\textbf{}]\label{T.G_i}
If $\Omega$ is an infinite set and $(G_i)_{i\in I}$ a chain
of subgroups {\rm(\!\<}or more generally, submonoids\/{\rm)}
of $S=\Sym(\Omega),$ with $\bigcup_{i\in I}\,G_i=
S$ and $|I|\leq|\<\Omega\<|,$ then $G_i=S$ for some $i\in I.$
\end{theorem}\begin{proof}
Lemma~\ref{L.oEfull} shows that $\Omega$ has a full moiety with respect
to some $G_i,$ hence assuming the $G_i$ are subgroups, Lemma~\ref{L.G&x}
shows that $S=\langle G_i\cup\{x\}\rangle$ for some $x\in S.$
Since the $G_j$ form a chain with union $S,$ there is some $j\geq i$
with $x\in G_j,$ hence $G_j\supseteq G_i\cup\{x\},$ hence $G_j=S.$

If the $G_i$ are merely submonoids, we apply the result of the
preceding paragraph to $(G_i\cap G_i^{-1})_{i\in I},$
which is clearly a chain of sub{\em groups} with union $S.$
\end{proof}

Now for our new result.

\begin{theorem}\label{T.U^n}
Suppose $\Omega$ is an infinite set, and $U$ a generating set for
$\Sym(\Omega)$ as a group {\rm(\!\<}respectively, as a monoid\/{\rm)}.
Then there exists a positive integer $n$
such that every element of $\Sym(\Omega)$ is represented by a group
word {\rm(\!\<}respectively, a monoid word\/{\rm)}
of length $\leq n$ in the elements of $U.$\qed
\end{theorem}\begin{proof}
Here it suffices to prove the monoid case, since the group words in the
elements of $U$ are just
the monoid words in the elements of $U\cup U^{-1}.$

So assume $U$ generates $S$ as a monoid.
For $i=1,2,\ldots\,,$
let $U_i=(U\cup\{1\})^i\cap(U^{-1}{\cup}\,\{1\})^i.$
By assumption the sets $(U\cup\{1\})^i$ have union $S,$ hence so do
their inverses, $(U^{-1}\cup\{1\})^i,$ hence so do the
intersections $U_i.$
Since $\aleph_0\leq|\<\Omega\<|,$ Lemma~\ref{L.oEfull} says that
$\Omega$ has a full moiety with respect to some $U_i.$
By Lemma~\ref{L.G&x} there exists $x\in S$ such that
$S=(U_i\<x)^7\<U_i ^2\<x\,\cup\,(x\<U_i )^7 x\<U_i^2,$ which we
see is contained in $(U_i\<\cup\<\{x\})^{17}.$
Taking a $j\geq i$ such that $x\in U_j,$ we get $(U\cup\{1\})^{17j}
=((U\cup\{1\})^j)^{17}\supseteq(U_i\cup\{x\})^{17}=S.$
\end{proof}

One may ask whether for a given set $\Omega,$ there is some single $n$
as in Theorem~\ref{T.U^n} that works for every generating set $U.$
To see that this is not so, let $\Omega=\Q/\Z,$ and let us
give this set the natural metric, of diameter $1/2,$ under which the
distance between two cosets of $\Z$ is the minimum of the
distances between their members, as real numbers.
(In other words, let us use the metric on $\Omega=\Q/\Z$ induced by the
arc-length metric on $\mathbb R/\Z.)$
Fixing an integer $n,$ let $U$ denote the set of permutations
of $\Omega$ which move all elements by distances $<1/(2n)$.
Clearly, $U^n\neq\Sym(\Omega).$
However, I claim that $U$ is a generating set for $\Sym(\Omega)$
as a group (and indeed, since it is closed under inverses, as a monoid).

Note first that if $\Sigma$
is the image in $\Omega$ of any interval of length $<1/(2n)$ in $\Q\<,$
then $U$ contains $S_{(\Omega-\Sigma)},$ the group of permutations
that act arbitrarily on $\Sigma$ and fix all elements outside it.
Now we can cover $\Omega$ with a finite number of successive
overlapping sets $\Sigma_1, \Sigma_2, \ldots, \Sigma_r$ of this sort,
and then use Lemma~\ref{L.S1S2} to conclude inductively that
$\langle\<U\rangle\supseteq
S_{(\Omega-(\Sigma_1\cup\ldots\cup\Sigma_i))}$
for $i=1,\ldots,r,$ hence that $\langle\<U\rangle\supseteq
S_{(\Omega-(\Sigma_1\cup\ldots\cup\Sigma_r))}=S_{(\varnothing)}=S,$
as claimed.
\vspace{6pt}

\section{Questions, examples, remarks,
and related literature.}\label{S.questions}

The statement of \cite[Theorem~1.1]{DM&PN} ${=}$~Theorem~\ref{T.G_i}
above, unlike that of Theorem~\ref{T.U^n}, depends on the cardinal
$|\<\Omega\<|.$
Let us, for the purposes of this section, weaken that statement to
one that holds independent of this cardinal, by using the obvious
lower bound $\aleph_0\leq|\<\Omega\<|.$
Then that theorem can be looked at as saying, for every infinite
set $\Omega,$ that
$\Sym(\Omega)$ belongs to the class of groups $G$ satisfying

\begin{xlist}\item\label{x.Gneqcup}
If $G$ is written as the union of a chain of subgroups
$G_0\leq G_1\leq\ldots$ indexed by $\omega,$ then for some $n,$ $G_n=G.$
\end{xlist}
Theorem~\ref{T.U^n}, similarly, says that
$\Sym(\Omega)$ belongs to the class of groups $G$ satisfying

\begin{xlist}\item\label{x.bd_n}
If $G$ is generated as a group by a subset $U,$ then for some
$n,$ every element of $G$ is represented by a group
word of length $\leq n$ in the elements of $U.$
\end{xlist}

Clearly~(\ref{x.Gneqcup}) also holds for all finitely generated groups,
and~(\ref{x.bd_n}) for all finite groups.
Thus, in a strange way, the groups $\Sym(\Omega)$
resemble finite groups.

Condition~(\ref{x.Gneqcup}) on a non-finitely-generated group is
commonly expressed by saying that $G$ has {\it uncountable cofinality,}
and is known to hold in many cases.
Some works on cofinalities of groups
are~\cite{PB}, \cite{MD+RG}, \cite{SK+JT}, \cite{GS}, \cite{STcof} and
\cite{ST_surv}; see also other papers cited in~\cite{MD+RG}.
Whether or not $G$ is finitely generated,~(\ref{x.Gneqcup})
is equivalent to the condition that the
fixed point set construction on $G$-sets commutes with direct limits
over countable index sets \cite[end of \S2]{dirlimfix}.
In fact, it was the wish to give in~\cite{dirlimfix} an example of a
non-finitely-generated group having this property that led me
to read~\cite{DM&PN}, eventually resulting in this note.

In response to a question posed in an earlier version of this note,
Droste and G\"obel~\cite{MD+RG} have obtained a general technique for
showing that certain sorts of structures containing
many isomorphic copies of themselves have automorphism groups
satisfying~(\ref{x.Gneqcup}) and~(\ref{x.bd_n});
in particular, they find that~(\ref{x.Gneqcup}) and~(\ref{x.bd_n})
hold for the groups of all self-homeomorphisms of
the Cantor set, of the rational numbers, and of the irrational numbers,
and for the group of Borel automorphisms of the real numbers.
Droste and Holland~\cite{MD+WCH} obtain the same conclusions for
the automorphism group of any doubly homogeneous totally ordered set,
and Tolstykh~\cite{VT}, \cite{VTrel} proves~(\ref{x.bd_n}),
and, insofar as it was not already known,~(\ref{x.Gneqcup}),
for the automorphism groups of infinite-dimensional vector spaces
and various sorts of relatively free groups.
In the final remark of~\cite{VT} he suggests an approach to showing
that the full automorphism groups of free objects in other well-behaved
varieties of algebras have this property.
Mesyan~\cite{ZM} obtains analogous results for endomorphism
rings of infinite direct sums and products of a module, and
Cornulier~\cite[Proposition~4.4]{YdC} does the same for
the Boolean ring of subsets of an infinite set;
he also shows in~\cite[Theorem~3.1]{YdC} that any
$\omega_1$\!-existentially closed group
satisfies~(\ref{x.Gneqcup}) and~(\ref{x.bd_n}).

Let us note some general facts about
conditions~(\ref{x.Gneqcup}) and~(\ref{x.bd_n}).
It is not hard to see that both properties are preserved under
taking homomorphic images, and that~(\ref{x.Gneqcup}) is also preserved
under group extensions and under passing to groups finitely
generated over $G.$
To get some similar results for~(\ref{x.bd_n}), we will want

\begin{lemma}\label{L.cosets}
Let $H<G$ be groups and $U$ a generating set for $G.$
Suppose that for some $n\geq 0,$ every right coset of $H$ in $G$
contains a group word of length $\leq n$ in the elements of $U.$
Then the set of elements of $H$ that can be written as words of length
$\leq 2n+1$ in the elements of $U$ generates $H.$
\end{lemma}\begin{proof}
Let $V$ be a set of right coset representatives for $H$ in $G$
consisting of words of length $\leq n$ in $U,$ with the coset $H$
represented by the element $1,$ and let $r:G\rightarrow V$
be the retraction collapsing each coset to its representative.
Let $W$ denote the set of elements of $H$ that can be written as words
of length $\leq 2n+1$ in the elements of $U.$

For any $v\in V$ and $u\in U\cup U^{-1},$ note that
$v\<u=(vu\;r(vu)^{-1})\<r(vu).$
Since $r(vu)$ by definition lies in the same right coset as $vu,$
the factor $vu\;r(vu)^{-1}$ lies in $H,$
and since $v$ and $r(vu),$ as members of $V,$ each have length
$\leq n$ in the elements of $U,$ that factor has length $\leq 2n+1,$
and so lies in $W.$
Thus, $V(\<U\cup U^{-1})\subseteq W\<V.$
It follows that $\bigcup_i W^i V$ is closed under right multiplication
by $U\cup U^{-1},$ hence equals the whole group $G.$
We now intersect the equation $\bigcup_i W^i V=G$ with $H.$
This has the effect of discarding elements on the
left-hand side having right factors from $V$
other than $1,$ and so gives $\bigcup_i W^i=H,$ completing the proof.
\end{proof}

We can now show that~(\ref{x.bd_n}) is preserved under group extensions.
Given a short exact sequence
$1\rightarrow H\rightarrow G\rightarrow E\rightarrow 1$ where $H$
and $E$ satisfy~(\ref{x.bd_n}), and a generating set $U$ for $G,$
the fact that $E$ satisfies~(\ref{x.bd_n})
yields an $n$ as in the hypothesis of Lemma~\ref{L.cosets}.
The conclusion of that lemma, combined with the
fact that $H$ satisfies~(\ref{x.bd_n}),
shows that all elements of $H$ can be written as
words of length $\leq m$ in the elements of $U$ for some $m.$
It follows that all elements
of $G$ can be written as words of length $\leq n+m.$
A similar application of that lemma shows that~(\ref{x.bd_n}) is
preserved under passing to overgroups in which $G$ has finite index.

Clearly, a countable group satisfies~(\ref{x.Gneqcup}) if and
only if it is finitely generated, while an infinite group
that is finitely generated can never satisfy~(\ref{x.bd_n}).
In particular,~(\ref{x.bd_n}) is not preserved under passing to groups
finitely generated over $G.$
Neither property is preserved under passing to normal subgroups,
since for $\Omega$ a countably infinite set, the subgroup of
$\Sym(\Omega)$ consisting of permutations that move only finitely
many elements is normal, but satisfies neither~(\ref{x.Gneqcup})
nor~(\ref{x.bd_n}).

In response to a question posed in an earlier version of
this preprint, Kh\'elif has announced in~\cite{AK}
that~(\ref{x.Gneqcup}), and also the conjunction
of~(\ref{x.Gneqcup}) and~(\ref{x.bd_n}), are preserved under passing
to subgroups of finite index, but that~(\ref{x.bd_n}) alone is not.

A question that, to my knowledge, is still unanswered is

\begin{question}\label{Q.ctbl&bd_n}
Are there any countably infinite groups that
satisfy~{\rm(\ref{x.bd_n})}?
\end{question}

One can show that every non-finitely-generated abelian
group can be mapped surjectively either onto a group $Z_{p^\infty},$ or
onto an infinite direct sum $Z_{p_1}\oplus Z_{p_2}\oplus...$\ for
(not necessarily distinct) prime numbers $p_1,\,p_2,\ldots\,.$
Neither of the latter sorts of groups satisfy~(\ref{x.Gneqcup})
or~(\ref{x.bd_n}), so no non-finitely-generated abelian group
has either of these properties.

It is shown in \cite{SK+JT} that the direct product of any family of
copies of a nonabelian finite simple group satisfies~(\ref{x.Gneqcup}).
In \cite[\S4]{YdC} this is strengthened to say that the direct product
of any family of copies of a finite perfect group
satisfies~(\ref{x.Gneqcup}) and~(\ref{x.bd_n}),
while~\cite[Lemma~3.5]{MD+RG} shows that the same is true of the
direct product of any family of copies of $\Sym(\omega)$.
These results suggest the general question of when~(\ref{x.Gneqcup})
and~(\ref{x.bd_n}) are inherited by products.
They are both inherited by finite products, by our observations
on group extensions.
That they are not always inherited by infinite products is shown by any
infinite direct product of nontrivial finite abelian groups:
the factors, as finite groups, satisfy both conditions, but the
product, a non-finitely-generated abelian group, satisfies neither.
The same technique can be used to show the failure of these conditions
for certain direct products of perfect groups, using the fact that an
infinite product of perfect groups may be non-perfect.
Kh\'elif~\cite{AK} has announced the existence of a direct product
of groups such that each factor satisfies~(\ref{x.bd_n}) and the
product is perfect but which still does not satisfy~(\ref{x.bd_n}).

We saw in the proof of Theorem~\ref{T.G_i} that~(\ref{x.Gneqcup})
implies the corresponding condition with ``subgroups'' replaced by
``submonoids''.
The same implication for the negations of these conditions is clear,
so the group and monoid versions of~(\ref{x.Gneqcup}) are equivalent.
For~(\ref{x.bd_n}), the proof of Theorem~\ref{T.U^n} showed
that the statement for generation as a monoid implied
the statement for generation as a group.
In this case, we can also get the converse {\em if} we
assume~(\ref{x.Gneqcup}).
For suppose the group $G$ satisfies~(\ref{x.Gneqcup})
and~(\ref{x.bd_n}), and that $U$ is a generating set for $G$
as a monoid.
As in the proof of Theorem~\ref{T.U^n},
let $U_i=(U\cup\{1\})^i\cap(U^{-1}{\cup}\,\{1\})^i$ $(i\in\omega).$
These sets form a chain with union $G,$ hence so do the
subgroups $\langle\<U_i\rangle,$ hence by~(\ref{x.Gneqcup}),
some $\langle\<U_j\rangle$ is equal to $G;$ so by~(\ref{x.bd_n}),
there is an $n$ such that all elements of $G$ are group words of length
$\leq n$ in elements of $U_j.$
By construction, $U_j$ is closed under inverses, so these group
words reduce to monoid words, so
$G=(U_j)^n\subseteq ((U\cup\{1\})^j)^n= (U\cup\{1\})^{jn},$ as claimed.
However, the results announced in~\cite{AK} on
subgroups of finite index imply that
there are groups satisfying~{\rm(\ref{x.bd_n})} but
not~{\rm(\ref{x.Gneqcup})}, leading to

\begin{question}\label{Q.combns}
Among groups satisfying~{\rm(\ref{x.bd_n})} but
not~{\rm(\ref{x.Gneqcup})}, are there any which do not satisfy the
analog of~{\rm(\ref{x.bd_n})} obtained by replacing ``generated as a
group'' and ``group word'' with ``generated as a monoid'' and
``monoid word''?
Any which {\em do} satisfy that condition?
\end{question}

We saw at the end of \S\ref{S.U^n} that though $\Sym(\Omega)$
satisfies~(\ref{x.bd_n}), there is no single $n$ that works
for all generating sets $U.$
On the other hand, Shelah~\cite{Sh} constructs an uncountable
group $G$ in which every generating set $U$ satisfies $U^{240}=G.$

We record two easily reformulations of the conjunction
of~{\rm(\ref{x.Gneqcup})} and~{\rm(\ref{x.bd_n})}.

\begin{lemma}\label{L.both}
For any group $G,$ the following conditions are equivalent:
\\[2pt]
{\rm(i)}~$G$ satisfies both~{\rm(\ref{x.Gneqcup})}
and~{\rm(\ref{x.bd_n})}.
\\[2pt]
{\rm(ii)}~If $U_0\subseteq U_1\subseteq\ldots\subseteq
U_i\subseteq\ldots$ is a chain of subsets of $G,$ indexed by the natural
numbers and having union $G,$ such that $U_i = (U_i)^{-1}$ for
all $i,$ and such that for all $i,j$ there exists some
$k$ with $U_i U_j\subseteq U_k,$ then some $U_i$ is equal to $G.$
\\[2pt]
{\rm(iii)}~If $L$ is a natural-number-valued function on $G$ such
that for all $g,h\in G$ one has $L(g^{-1})=L(g)$
and $L(gh)\leq L(g)+L(h),$ then $L$ is bounded above.\qed
\end{lemma}

Functions $L$ as in~(iii) above are studied in geometric group
theory.
(Further conditions are generally assumed,
in particular, that the number of $g\in G$ with
$L(g)\leq n$ is finite and grows at most exponentially in $n.$
Under this assumption, $L$ is known to be approximable
by the function giving the length of $g$ with respect to
some finite generating set of an overgroup of $G$~\cite{Olsh}.
Of course, such a finiteness condition cannot be satisfied when
$G$ is an uncountable group such as we are considering.)
In \cite[\S2]{YdC}, condition~(iii) above is translated
as saying that every isometric action of $G$ on a metric space has
bounded orbits, and it is deduced that isometric actions of~$G$
on certain sorts of metric spaces must in fact have fixed points.

\vspace{6pt}
A property of finite groups $G$ that does {\em not} hold for
the groups $\Sym(\Omega)$ is

\begin{xlist}\item\label{x.genasmnd}
Every subset of $U\subseteq G$ which generates $G$ as a group generates
$G$ as a monoid.
\end{xlist}
To see the failure of~(\ref{x.genasmnd}) in $\Sym(\Omega)$ for countable
$\Omega,$ we may take for $U$ a submonoid $M$ as in the next result.

\begin{lemma}\label{L.updwnup}
Let $\Omega$ be any countable totally ordered set without least or
greatest element, and let $M$ be the monoid $\{g\in\Sym(\Omega):
(\<\forall\<\alpha\in\Omega)\,\,\alpha\<g\geq\alpha\}.$
Then $\Sym(\Omega)=M\<M^{-1}\!\<M.$
\end{lemma}\begin{proof}
Given $f\in\Sym(\Omega),$ I claim that we can find
$g,h\in\Sym(\Omega)$ such that for all $\alpha\in\Omega,$

\begin{xlist}\item\label{x.updwnup}
$\alpha\,\leq\,\alpha\<g\,\geq\,\alpha\<h\,\leq\,\alpha f.$
\end{xlist}
From this it will follow that $g,\ h^{-1}g,\ h^{-1}f\in M,$
whence $f=g\<(h^{-1}g)^{-1}(h^{-1}f)\in M\<M^{-1}\!\<M,$ as desired.

To get $g$ and $h$ satisfying~(\ref{x.updwnup}),
let an enumeration of the elements of $\Omega$ be chosen.
(We will not introduce a notation for this enumeration, but simply
speak of ``the first element with respect to our enumeration
such that ...''.
In particular, ``$\leq$'' will continue to denote the given ordering of
$\Omega,$ not the ordering corresponding to our enumeration.)
We shall now construct recursively for each $i\in\omega$ a $3$-tuple
$(\alpha_i,\,\beta_i,\,\gamma_i)$ of elements
of $\Omega;$ these will eventually be the $3$-tuples
$(\alpha,\,\alpha\<g,\,\alpha\<h)$ $(\alpha\in\Omega).$

To define $(\alpha_i,\,\beta_i,\,\gamma_i),$ let us,
if $i$ is divisible by $3,$ take for $\alpha_i$
the first element of $\Omega$ with respect to our enumeration which
has not been chosen as $\alpha_j$ at any previous step (i.e.,
for $0\leq j<i),$ let us then take for $\beta_i$ any element of $\Omega$
which was not chosen as $\beta_j$ at a previous step and which
is $\geq\alpha_i,$ and for $\gamma_i$
any element not chosen as $\gamma_j$ at a previous step which
is both $\leq\beta_i$ and $\leq\alpha_i\<f.$
On the other hand, if $i\equiv 1~(\mathrm{mod}~3),$ we start
by taking for $\beta_i$ the first element of $\Omega$ with respect
to our enumeration which was not chosen as $\beta_j$ at any
previous step, then take for $\alpha_i$ any element not previously
chosen as an $\alpha_j$ which is $\leq\beta_i,$ and for $\gamma_i$
any element not previously chosen as a $\gamma_j$ which
is both $\leq\beta_i$ and $\leq\alpha_i\<f.$
Finally, when $i\equiv 2~(\mathrm{mod}~3),$ we take for
$\gamma_i$ the first element with respect to our enumeration
that has not yet been used as a $\gamma_j,$ for $\alpha_i$
any element not previously chosen as an $\alpha_j$ such
that $\alpha_i\<f\geq\gamma_i$ (which is possible because there are
infinitely many elements $\geq\gamma_i),$ and for $\beta_i$ any
element not previously chosen in that role which
is both $\geq\alpha_i$ and $\geq\gamma_i.$

Clearly this construction uses each element of $\Omega$ once and only
once in each position; hence the set of pairs $(\alpha_i,\,\beta_i)$ is
the graph of a permutation $g\in\Sym(\Omega),$ and the set of pairs
$(\alpha_i,\,\gamma_i)$ is the graph of a
permutation $h\in\Sym(\Omega).$
The conditions we imposed on our choices at each step insure
that~(\ref{x.updwnup}) holds.
\end{proof}

(In the above lemma we can drop the countability assumption, if
we replace the hypothesis of no least or greatest element by
the assumption that each element of $\Omega$ has
$|\<\Omega\<|$ elements above it and $|\<\Omega\<|$ elements below it.
The assumption that the ordering is total can also be
weakened to say that it is upward and downward directed.)

It would be interesting to know what sorts of groups
satisfy~(\ref{x.genasmnd}), other than those whose elements all
have finite order.
One such group is the infinite dihedral group.

My first attempts to prove the statements about sub{\em monoids} of
$S$ in Theorems~\ref{T.G_i} and~\ref{T.U^n} revolved
around trying to remove from Lemmas~\ref{L.S1S2}
and~\ref{L.G&x} the hypotheses on closure under inverses.
Those hypotheses are required by our proof of Lemma~\ref{L.S1S2}, since
inverses are needed to form commutators and conjugates.
Though a different method eventually gave the monoid cases of
those theorems, it would still be interesting to know the answer to

\begin{question}\label{Q.mS1S1}
Are versions of Lemma~\ref{L.S1S2} and/or~\ref{L.G&x}
{\rm(\!\<}possibly using longer products of $x,$ $U$ and $V)$
true without the hypothesis that $U$ and $V$ are closed under inverses?
\end{question}

We shall note a
weak result in this direction at the end of the next section.

We remark, finally, that it is easy to adapt
the method of proof of Theorem~\ref{T.G_i} to give an
apparently more general statement, in which the {\em chain} of subgroups
$G_i$ indexed by a set $I$ of cardinality $\leq |\<\Omega\<|$
is replaced by any {\em directed system} of subgroups $G_i$ indexed
by such a set, again having $S$ as union.
However, that result in fact follows easily from the theorem as stated.
For given such a directed system, let $\{G_i:i\in\kappa\}$ be a
subset thereof whose union generates $S,$ having least cardinality
$\kappa$ among such subsets, and indexed by that cardinal.
Then $\{\langle\,\bigcup_{i<j} G_i\rangle:j<\kappa\}$ forms a chain of
{\em proper} subgroups of $S,$ which, unless $\kappa$ is finite,
will have union $S.$
This would contradict Theorem~\ref{T.G_i}; so $\kappa$ is finite,
and the finite family $\{G_i:i\in\kappa\}$ will be majorized by a member
of the original directed system, which thus equals $S.$

\vspace{6pt}
Further results on the groups $\Sym(\Omega)$ are obtained
in~\cite{Sym_Omega:2}.

\section{Appendix:  Writing every element of $\Sym(\Omega)$ as a
commutator.}\label{S.rplt}
In proving Lemma~\ref{L.S1S2}, we called on the result of~\cite{OO}
that every element of an infinite symmetric group is a commutator.
Now a commutator is an element obtained by dividing an element
by a conjugate, $(g^{-1}hg)^{-1}h;$ moreover, it is clear that in
a symmetric group, every element is conjugate to its inverse; so
a commutator in a symmetric group can be described as a product of
two elements in the same conjugacy class.
Almost a decade after~\cite{OO} appeared, it was
shown that in an infinite symmetric group there in
fact exist single conjugacy classes whose square is the whole group.
In a series of papers by several authors, culminating
in~\cite{GM} (which describes this history), the conjugacy
classes having this property were precisely characterized.

We shall give a self-contained proof of this property for
one such conjugacy class in Lemma~\ref{L.rplt.rplt} below, then
use the fact that $\Sym(\Omega)$ is the square of a conjugacy class
to get a result related to Question~\ref{Q.mS1S1}.

\begin{definition}\label{D.replete}
For $\Omega$ an infinite set, we shall call an element
$f\in\Sym(\Omega)$ {\em replete} if it has $|\<\Omega\<|$
orbits of each positive cardinality $\leq\aleph_0$
{\rm(}including $1).$
For a subset
$\Sigma\subseteq\Omega$ of cardinality $|\<\Omega\<|,$ we shall
say that $f$ is replete on $\Sigma$ if $\Sigma f=\Sigma$ and the
restriction of $f$ to $\Sigma$ is a replete permutation of $\Sigma.$
\end{definition}

Note that a permutation of $\Omega$ that is replete on a subset
$\Sigma\subseteq\Omega$ of cardinality $|\<\Omega\<|$ is necessarily
replete on~$\Omega\<.$

The replete permutations of $\Omega$ clearly form a conjugacy class,
so it will suffice to prove

\begin{lemma}\label{L.rplt.rplt}
Every permutation $f$ of an infinite set $\Omega$ is the product
of two replete permutations.
\end{lemma}\begin{proof}
Given $f,$ let us choose a moiety $\Sigma_0$ of $\Omega$
such that $f$ moves only finitely many elements
from $\Sigma_0$ to $\Omega-\Sigma_0$ or from $\Omega-\Sigma_0$ to
$\Sigma_0.$
That there exists such a $\Sigma_0$ is immediate
if $\Omega$ is uncountable, for
we can break $\Omega$ into two families of $|\<\Omega\<|$ orbits each,
and take for $\Sigma_0$ the union of one of these families.
If $\Omega$ is countable, we can use the same method if $f$ has
infinitely many orbits, and can also get the same conclusion in an
obvious way if $f$ has more than one infinite orbit.
If it has exactly one infinite orbit, $\alpha_0\<\langle f\<\rangle,$
and finitely many finite orbits, then we can take
$\Sigma_0=\{\alpha_0 f^n:n\leq 0\};$ clearly $f$ moves exactly
one element out of $\Sigma_0,$ and none into it.

After choosing $\Sigma_0,$
let us split $\Omega-\Sigma_0$ into two disjoint moieties
$\Sigma_1$ and $\Sigma_2,$ so that $\Sigma_1$ contains the finitely
many elements of $(\Sigma_0 f\,\cup\,\Sigma_0 f^{-1})-\Sigma_0.$
I claim that if $g_0$ is any permutation of $\Sigma_0$ and
$g_2$ any permutation of $\Sigma_2,$
then there exists a permutation $h$ of $\Omega$
such that $fh$ agrees on $\Sigma_0$ with $g_0,$ while
$h$ agrees on $\Sigma_2$ with $g_2.$
Indeed, this pair of conditions specifies the values of $h$ on the
two disjoint sets $\Sigma_0 f$ and $\Sigma_2$ in a one-to-one
fashion, and both the set on which it leaves $h$ unspecified
and the set of elements that it does not specify as values
for $h$ are of cardinality $|\<\Omega\<|.$
Hence the former set can be mapped bijectively to the latter, and the
resulting bijection will serve to complete the definition of $h.$
(We have used $\Sigma_1$ as a ``Hilbert's Hotel'' in
case $\Sigma_0 f\neq\Sigma_0.)$

Now if we take for $g_0$ and $g_2$ replete permutations
of $\Sigma_0$ and $\Sigma_2$ respectively,
then $h$ will be replete on $\Sigma_2,$ hence replete,
while $fh$ will be replete on $\Sigma_0,$ hence replete.
Thus $f=(fh)\<h^{-1}$ is a product of two replete permutations, as
we wished to show.
\end{proof}

(In~\cite[Theorem~3.1(a)]{MD85} the same result is proved
for a different
conjugacy class, that of permutations with $|\<\Omega\<|$ infinite
orbits and no finite orbits, also by an unexpectedly simple argument.
However, it takes some work to isolate that argument from the
lengthier proofs of other results that are being carried out
there simultaneously.
A different sort of generalization of Ore's result is obtained
in~\cite{D+M} and papers cited there, which characterize the group
words $w$ which are ``universal'' in infinite permutation groups, in
the sense established for the word $x^{-1} y^{-1}x\<y$ by Ore's result.)
\vspace{6pt}

We shall now show that the proofs of Lemmas~\ref{L.S1S2} and~\ref{L.G&x}
can be adapted to the situation where $U$ and $V$ are not assumed
closed under inverses if we allow ourselves to use,
along with multiplication, the right conjugation operation
\begin{xlist}\item\label{x.conj}
$g^h\;=\;h^{-1}g\,h\<.$
\end{xlist}

\begin{lemma}[{\rm cf.\ Lemma~\ref{L.G&x}}\textbf{}]\label{L.G&x&y}
Let $\Omega$ be an infinite set, and
$U\subseteq S=\Sym(\Omega)$ a subset with respect to which
some moiety of $\Omega$ is full.
Then there exist $x,y\in S,$ with $x$ of order $2,$ such that

\begin{xlist}\item\label{x.UxUxUy^Uy^U}
$S\;=\;(Ux)^3(y^U)^2\<x\;\cup\;x\<(Ux)^3(y^U)^2.$
\end{xlist}
\end{lemma}
\noindent
{\it Sketch of proof.}
Let $\Sigma_1$ be a full moiety for $U,$ let $\Sigma_2$ be
a moiety such that $\Sigma_1\cap\Sigma_2$ is a moiety and
$\Sigma_1\cup\Sigma_2=\Omega,$
let $x\in S$ be an involution such that $\Sigma_1\,x=\Sigma_2,$ and
let $y\in S_{(\Sigma_2)}$ be an element such that the
group $S_{(\Sigma_2)}\cong\Sym(\Omega-\Sigma_2)$ is the square of
the conjugacy class of $y$ in that group.
Such a $y$ exists by
Lemma~\ref{L.rplt.rplt} above, or the results in the papers cited.
Combining this property of $y$ with the fact that $\Sigma_1,$
and hence its subset $\Omega-\Sigma_2,$ is full with respect
to $U,$ we conclude that
$S_{(\Sigma_2)}\subseteq (y^U)^2;$ hence, conjugating by $x,$
\begin{xlist}\item\label{x.y^U^2}
$S_{(\Sigma_1)}\subseteq x\<(y^U)^2\<x.$
\end{xlist}

Now let $V=x\,Ux.$
Since $\Sigma_1$ is a full moiety under $U,$ $\Sigma_2$ will
be a full moiety under $V.$
Using the technique of proof of Lemma~\ref{L.S1S2}, with~(\ref{x.y^U^2})
in place of~(\ref{x.VUVUVV}), we get~(\ref{x.UxUxUy^Uy^U}).
\qed\vspace{6pt}

I am grateful to Peter Biryukov, Manfred Droste, Dugald Macpherson,
Gadi Moran, Peter Neumann and Marc Rieffel for helpful
information on the relevant literature.

George M. Bergman\\
Department of Mathematics\\
University of California\\
Berkeley, CA 94720-3840\\
USA\\
\bigskip
gbergman@math.berkeley.edu
\end{document}